\documentclass[preprint,12pt]{elsarticle}

\usepackage{amsmath,amssymb,amsthm}
\usepackage{graphicx}
\usepackage{tikz}
\usetikzlibrary{arrows.meta,calc}
\usepackage{siunitx}
\usepackage{bm}
\usepackage{lineno}         

\newtheorem{theorem}{Theorem}[section]

\newtheorem{corollary}  [theorem]{Corollary}

\theoremstyle{definition}
\newtheorem{remark}[theorem]{Remark}

\newcommand{\vect}[1]{\bm{#1}}

\newcommand{\grad}{\nabla}
\newcommand{\Div}{\nabla\!\cdot}
\newcommand{\curl}{\nabla\times}

\journal{Applied Mathematics and Computation}

\begin{document}

\begin{frontmatter}

\title{Algebraic Reductibility Experiments of RANS-Inspired Equations}

\author[ibero]{Carla Valencia\corref{cor1}}
\ead{carla.valencia@ibero.mx}

\author[ibero]{Sebastián Velasco}
\author[ibero,vera]{Manuel Romero de Terreros}

\cortext[cor1]{Corresponding author}

\affiliation[ibero]{%
  organization={Universidad Iberoamericana Ciudad de México},
  city={Mexico City},
  postcode={01219},
  country={Mexico}}

\affiliation[vera]{%
  organization={Vera Strata Research},
  city={Mexico City},
  postcode={05100},
  country={Mexico}}

\begin{abstract}
Prior to any statistical averaging we derive a rotational form of the
Reynolds-Averaged Navier–Stokes (RANS) equations, which eliminates the
pressure and exposes a velocity–vorticity interplay governed by
\[
\partial_t(\boldsymbol{\omega}+\tilde{\boldsymbol{\omega}})
      +(\vect{v}\!\cdot\!\nabla)\boldsymbol{\omega}
      +(\tilde{\vect{v}}\!\cdot\!\nabla)\tilde{\boldsymbol{\omega}}
      +(\vect{v}\!\cdot\!\nabla)\tilde{\boldsymbol{\omega}}
      +(\tilde{\vect{v}}\!\cdot\!\nabla)\boldsymbol{\omega}
      -\nu\Delta(\boldsymbol{\omega}+\tilde{\boldsymbol{\omega}})=\vect{0}.
\]
All terms are differential polynomials, so the system generates a
differential–algebraic ideal.  Using the \emph{Rosenfeld–Gröbner} algorithm
we obtain an equivalent triangular hierarchy whose first equation
involves a single variable, the second at most two, and so on.  This
decoupling clarifies how prescribed mean-flow data drive the turbulent
fluctuations and provides a hierarchy-ready foundation for
physics-informed or physics-embedded neural networks.  Energy estimates
in Sobolev spaces complement the algebraic reduction and establish
local well-posedness when the initial kinetic energy of the velocity
and its curl is finite.  The joint algebraic–energetic framework thus
offers a pressure-free, computationally economical platform for
data-driven turbulence analysis.
\end{abstract}

\begin{graphicalabstract}
  \centering
  \includegraphics[width=0.8\textwidth]{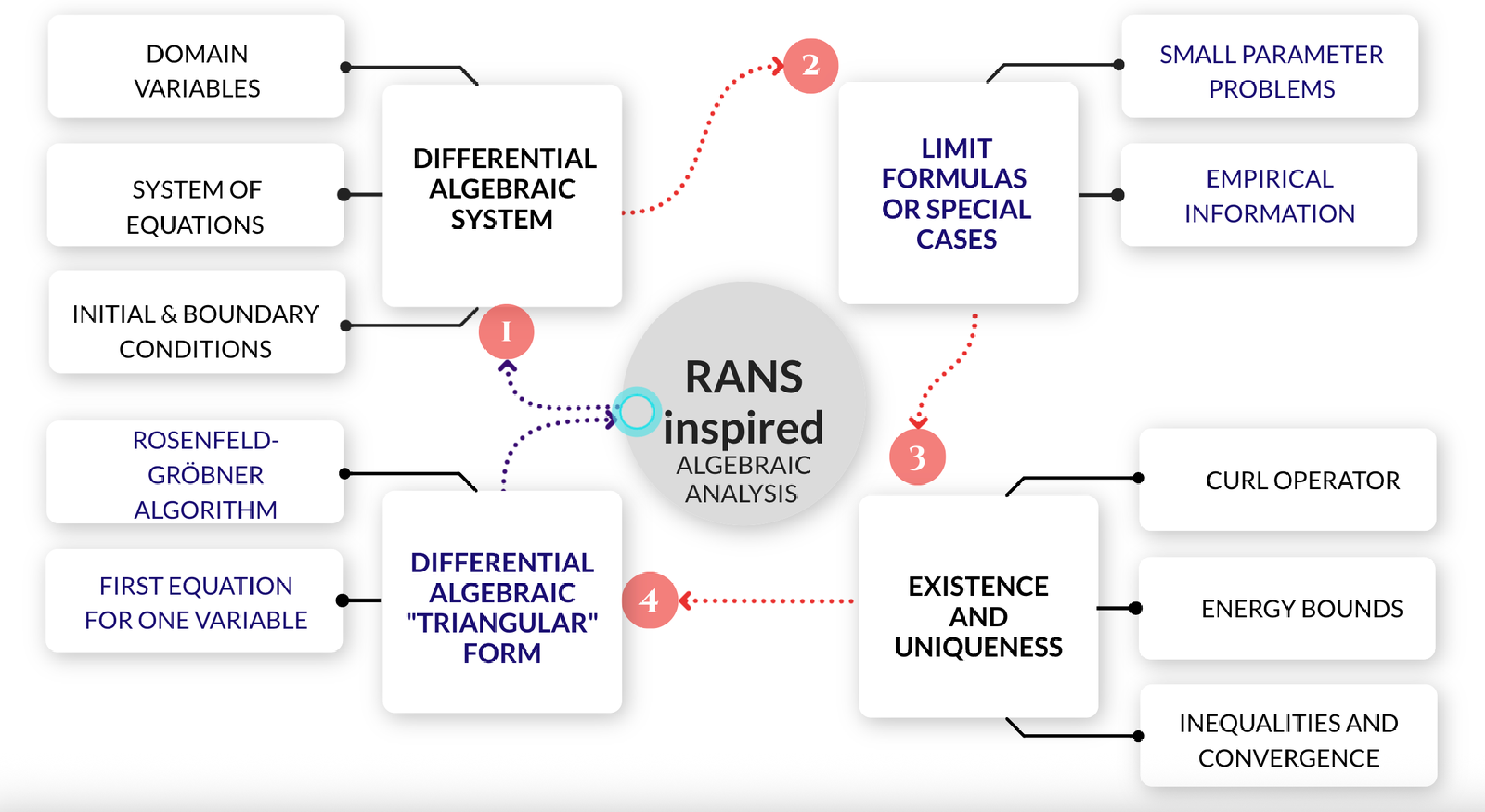}
\end{graphicalabstract}

\begin{highlights}
\item Pressure‐free rotational formulation
\item \emph{Rosenfeld–Gröbner} yields a triangular hierarchy
\item Hierarchy accelerates training of physics-informed neural networks
\end{highlights}

\begin{keyword}
Turbulence modelling \sep Navier–Stokes equations \sep
Differential algebra \sep Rosenfeld–Gröbner \sep
Physics-informed neural networks
\\[2pt]
\MSC[2020] 35Q30 \sep 76F55 \sep 12H05 \sep 35A01 \sep 68T07
\end{keyword}

\end{frontmatter}

\section{Introduction}\label{sec:intro}
The study of turbulent flow, is pivotal 
in engineering applications,
and the RANS equations remain central in its study \cite{Tominaga2024,Meneveau2019}.
Turbulent motions in natural and built environments are rarely
homogeneous: surface roughness patches, temperature gradients, and
obstacles introduce spatial inhomogeneities that deflect a
mean flow and seed multi-scale vortical structures.
Capturing their dynamics is crucial for applications that
range from atmospheric boundary-layer modelling to
fluid–structure-interaction (FSI) assessments in extreme-event
engineering.

Throughout this work we assume that the fluctuation field
$\tilde{\vect{v}}$ possesses a clear physical interpretation:
It captures the velocity deviations generated by environmental
inhomogeneities, such as surface roughness, mosaics, temperature gradients and obstacles.
arrays, etc.).  
One, $\vect{v}$ corresponds to the classical velocity, 
and the other, $\tilde{\vect{v}}$, 
to the deflection caused by the active grid, both with null boundary conditions. 
This viewpoint  provides a family of models able to 
describe and explain recent active-grid wind-tunnel
experiments, in which controlled boundary perturbations produce
atmospheric-like turbulence at unprecedented Reynolds numbers
\cite{Neuhaus2020,Neuhaus2021, Kildal2023, Kistner2024, Angriman2024}.

\subsection*{1.1  A pressure-free rotational RANS type expression}

Starting from the instantaneous Navier–Stokes system we perform a
two-scale split \(\vect{u}=\vect{v}+\tilde{\vect{v}}\) before
any averaging and apply the curl operator.  
The pressure term disappears, yielding the vorticity
balance:
\begin{equation}
  \partial_t(\boldsymbol{\omega}+\tilde{\boldsymbol{\omega}})
  +(\vect{v}\!\cdot\!\grad)\boldsymbol{\omega}
  +(\tilde{\vect{v}}\!\cdot\!\grad)\tilde{\boldsymbol{\omega}}
  +(\vect{v}\!\cdot\!\grad)\tilde{\boldsymbol{\omega}}
  +(\tilde{\vect{v}}\!\cdot\!\grad)\boldsymbol{\omega}
  -\nu\nabla^{2}(\boldsymbol{\omega}+\tilde{\boldsymbol{\omega}})
  =\mathbf{0},
  \label{eq:rotrot}
\end{equation}
where \(\boldsymbol{\omega}=\curl\vect{v}\) and
\(\tilde{\boldsymbol{\omega}}=\curl\tilde{\vect{v}}\).
Equation \eqref{eq:rotrot} features one less variable than its
primitive-variable counterpart and displays, term-by-term, how the
mean profile and its curl steer the fluctuating vorticity and
vice-versa through distinct cross-convective actions.

\subsection*{1.2  Differential-algebraic hierarchy}

All operators in~\eqref{eq:rotrot} act polynomially on the dependent
variables, so the system forms a differential-algebraic ideal.
Implementing the Rosenfeld–Gröbner algorithm
\cite{Boulier1994,Boulier1995,BoulierLemaire2000,Boulier2009} we generate an
equivalent triangular system
\[
\mathcal{T}_1(\psi_1)=0,\quad
\mathcal{T}_2(\psi_1,\psi_2)=0,\quad
\mathcal{T}_3(\psi_1,\psi_2,\psi_3)=0,\;\ldots
\]
such that each block \(\mathcal{T}_k\) depends on at most \(k\)
unknowns \(\{\psi_1,\dots,\psi_k\}\).  
This cascading structure is particularly attractive for
physics-informed neural networks (PINNs) or
physics-embedded neural networks (PENNs), which can train each block
separately, propagate automatic-differentiation residuals efficiently,
and enforce hierarchical consistency.

\subsection*{1.3  Outline and contributions}

\noindent
\textbf{Pressure elimination.}  
      Curling the momentum equations yields \eqref{eq:rotrot},
      reducing the state dimension and sharpening the velocity-only
      perspective.
      \vspace{3mm}

\noindent      
\textbf{Transparent coupling.}  
      The cross-advective terms in \eqref{eq:rotrot} make explicit
      how the mean field’s initial/boundary data drive the fluctuating
      vorticity and vice-versa.
      \vspace{3mm}
      
\noindent
\textbf{Triangular reduction.}  
      A Rosenfeld–Gröbner procedure produces a hierarchical 
      representation with progressively fewer coupled variables,
      streamlining future neural-network solvers.
\vspace{3mm}

Because every term is a differential polynomial, the system defines a
differential ideal.  Section~\ref{sec:alg-reduction} applies the
Rosenfeld–Gröbner algorithm
\cite{BoulierLemaire2000,Boulier2009} to obtain equivalent
triangular forms when for suitable cases.  Section~\ref{sec:energy} derives Sobolev energy
identities that combine with the algebraic reduction to give local
existence and uniqueness in two and three dimensions.  Conclusions and
future directions close the paper.

\section{Problem statement}\label{sec:problem}

Let \(\Omega=(0,L_1)\times(0,L_2)\times(0,L_3)\subset\mathbb{R}^3\)
with \(L_i>0, \ i=1,2,3\).  We split the velocity field as
\(\bm{u}=\bm{v}+\tilde{\bm{v}}\) with
\(\Div\bm{v}=\Div\tilde{\bm{v}}=0\) and homogeneous Dirichlet data
\(\bm{v}=\tilde{\bm{v}}=\vect{0}\) on \(\partial\Omega\).
Define
\(\boldsymbol{\omega}:=\curl\bm{v}\) and
\(\tilde{\boldsymbol{\omega}}:=\curl\tilde{\bm{v}}\).

\begin{remark} The following results are stated so that they may be extended to \emph{starshaped} 
real domains $\Omega$ (where any ray with an origin in a fixed point has a unique common point with
the topological boundary of the domain)
and unknown solutions in a Sobolev space $W^k_p\left( \Omega \right)$, 
with $p\in \mathbb{N}-\{0\}$. This is possible thanks to the fact that
the space of infinitely differentiable functions defined on the topological closure 
$\bar{\Omega}=\Omega \cup \partial \Omega$ of a starshaped domain $\Omega$,
denoted as $C^{\infty}(\bar{\Omega})$, is dense in $W^k_p\left(\Omega\right)$,
where $\partial \Omega$ denotes the topological boundary of $\Omega$.
We denote $H^{m}(\Omega)
   \;:=\;
   W_2^m(\Omega)$.
\end{remark}

\begin{remark}
While a necessary and sufficient condition for a distribution in a Sobolev space 
to be the limit of a sequence of analytic functions is unknown 
it is possible to assume that the domain is starshaped
with respect to a point
to say that the space $C^{\infty}(\bar{\Omega})$  on a domain $\Omega$ 
is dense in the Sobolev space $W^k_p(\Omega)$ of distributions on $\Omega$ with derivatives 
of order $0\leq l \leq k$ in the Lebesgue space $L_p(\Omega)$, for $p\in[1,\infty)$.  
\end{remark}

\begin{center}
\begin{tikzpicture}[line join=round,
                    axis/.style   ={dashed, thick, red!70},
                    dashedblk/.style={dashed, thick},
                    scale=1.2]%
  \tikzset{x={(0.9cm,0.35cm)}, y={(-0.9cm,0.35cm)}, z={(0cm,1cm)}}

  \def\Lx{3}  
  \def\Ly{2}  
  \def\Lz{1.5}

  \coordinate(O)  at (0,0,0);
  \coordinate(A)  at (\Lx,0,0);
  \coordinate(B)  at (\Lx,\Ly,0);
  \coordinate(C)  at (0,\Ly,0);
  \coordinate(D)  at (0,0,\Lz);
  \coordinate(E)  at (\Lx,0,\Lz);
  \coordinate(F)  at (\Lx,\Ly,\Lz);
  \coordinate(G)  at (0,\Ly,\Lz);
  \coordinate(H)  at (0,0,1.7) ;
  \draw[thick] (O)--(A)--(B)--(C)--cycle;          
  \draw[thick] (D)--(E)--(F)--(G)--cycle;          
  \draw[thick] (O)--(D) (A)--(E) (B)--(F) (C)--(G);

  \path ($(O)!0.55!(F)$) node[scale=1.4]{$\Omega$};

  \draw[dashedblk] (O)--(B) (O)--(C) (O)--(D);

  \draw[axis,->] (O) --+ (0,0.8,0)   node[left] {$x_{1}$};
  \draw[axis,->] (O) --+ (0.8,0,0)    node[right]{$x_{2}$};
  \draw[axis,->] (O) --+ (0,0,1.0)    node[above]{$x_{3}$};

  \node[below left]  at (C) {$L_{1}$};
  \node[right]       at (A) {$L_{2}$};
  \node[above left]  at (H) {$L_{3}$};

\end{tikzpicture}
\end{center}

\begin{remark}
We may assume that the domain is a rectangular prism $\Omega$
and $\vect{v}$ is an element of the space $C^{\infty}(\bar{\Omega})$  
that satisfies both the Leibinz rule for the product derivative, 
the Gauss-Green formula \cite[p. 699]{Piskunov1969}, 
and the generalisation of Schwarz's Theorem \cite[p. 280]{Piskunov1969} 
to exchange the order of the generalised second partial derivatives 
with respect to the spatial and temporal independent variables \cite{Mazya1997,Mazya2013},
and perform an integration by parts process when needed. 
All differential identities below use that  
\(
  \vect{v},\tilde{\vect{v}}\in C^\infty(\overline{\Omega})
\),
vanish on \(\partial \Omega\) and are divergence-free;
hence surface integrals produced by integration-by-parts cancel,
and Schwarz’ theorem lets us swap mixed derivatives.
\end{remark}


\begin{theorem}[Rotational two–scale Navier–Stokes identity]
\label{thm:rot-two-scale}
Fix a rectangular box  
\[
   \Omega=(0,L_{1})\times(0,L_{2})\times(0,L_{3})\subset\mathbb{R}^{3},
   \qquad L_{1},L_{2},L_{3}>0,
\]
and a finite time horizon \(0<T<\infty\).
Let the pressure \(p\) and the velocity fields
\(
  \bm{v},\;\tilde{\bm{v}}\in
  C^{1}\!\bigl([0,T);C^{\infty}(\overline{\Omega};\mathbb{R}^{3})\bigr)
\)
satisfy, for every \(t\in[0,T)\),
\begin{enumerate}[\upshape(i)]
\item \textbf{Boundary conditions\,:}  
      \(\bm{v}=\tilde{\bm{v}}=\bm{0}\) on \(\partial\Omega\);
\item \textbf{Incompressibility\,:}  
      \(\nabla\!\cdot\!\bm{v}=0\) and \(\nabla\!\cdot\!\tilde{\bm{v}}=0\) in \(\Omega\);
\item \textbf{Momentum balance\,:}
      \begin{equation}\label{ns}
        \partial_{t}\!\left(\bm{v}+\tilde{\bm{v}}\right)
        +\bigl[(\bm{v}+\tilde{\bm{v}})\!\cdot\!\nabla\bigr]
          (\bm{v}+\tilde{\bm{v}})
        \;=\;
        \nu\,\Delta\!\left(\bm{v}+\tilde{\bm{v}}\right)-\nabla p
        \quad\text{in }\Omega,
      \end{equation}
      where \(\nu>0\) is the kinematic viscosity.
\end{enumerate}

\smallskip
Define the mean and fluctuating vorticities
\[
   \boldsymbol{\omega}:=\nabla\times\bm{v},
   \qquad
   \tilde{\boldsymbol{\omega}}:=\nabla\times\tilde{\bm{v}} .
\]
Then for all \((\bm{x},t)\in\Omega\times(0,T)\) the pair
\(\bigl(\boldsymbol{\omega},\tilde{\boldsymbol{\omega}}\bigr)\)
satisfies the \emph{rotational two–scale equation}
\begin{equation}
\label{eq:rot-two-scale}
\partial_{t}\bigl(\boldsymbol{\omega}+\tilde{\boldsymbol{\omega}}\bigr)
+\underbrace{(\bm{v}\!\cdot\!\nabla)\boldsymbol{\omega}
           +(\tilde{\bm{v}}\!\cdot\!\nabla)\tilde{\boldsymbol{\omega}}
           +(\bm{v}\!\cdot\!\nabla)\tilde{\boldsymbol{\omega}}
           +(\tilde{\bm{v}}\!\cdot\!\nabla)\boldsymbol{\omega}}_{\text{four distinct advection terms}}
-\nu\,\Delta\bigl(\boldsymbol{\omega}+\tilde{\boldsymbol{\omega}}\bigr)
=\bm{0}.
\end{equation}
In particular, the pressure has been eliminated and the dynamics are
expressed solely in terms of the velocity fields and their curls.
\end{theorem}

\begin{proof}
Apply \(\nabla\times(\,\cdot\,)\) to the momentum balance in
\eqref{ns} and use \(\nabla\times\nabla p=\boldsymbol{0}\):
\[
   \partial_{t}(\boldsymbol{\omega}+\tilde{\boldsymbol{\omega}})
   \;+\;
   \nabla\times\!\bigl((\vect{v}+\tilde{\vect{v}})\!\cdot\!\nabla\bigr)
        (\vect{v}+\tilde{\vect{v}})
   \;=\;
   \nu\Delta(\boldsymbol{\omega}+\tilde{\boldsymbol{\omega}}).
\]

\emph{Convection term.}
For any smooth vector fields \(\boldsymbol{a},\boldsymbol{b}\) with
\(\nabla\!\cdot\!\boldsymbol{a}=0\),
\(
   \nabla\times\bigl((\boldsymbol{a}\!\cdot\!\nabla)\boldsymbol{b}\bigr)
   =(\boldsymbol{a}\!\cdot\!\nabla)(\nabla\times\boldsymbol{b})
\).
Splitting \(\vect{v}+\tilde{\vect{v}}\) yields the four advection
terms visible in \eqref{eq:rot-two-scale}.

\emph{Boundary-term elimination.}
Integrating each component over \(R\) and integrating by parts,
all surface contributions vanish because
\(\vect{v}=\tilde{\vect{v}}=\boldsymbol{0}\) on \(\partial R\).
Hence the integral identity
\[
   \int_{\Omega}\bigl(\text{LHS of \eqref{eq:rot-two-scale}}\bigr)\,\mathrm d\boldsymbol x
   =\boldsymbol{0}
\]
holds.  
By the mean-value theorem for integrals there exists
\(\boldsymbol{x}_{0}\in \Omega \) where the integrand itself vanishes;
smoothness then forces equality to hold pointwise, yielding
\eqref{eq:rot-two-scale}.
\end{proof}

\begin{remark}
Equation~\eqref{eq:rot-two-scale} is obtained under the \emph{homogeneous
(no–slip)} boundary conditions
\[
  \vect{v}\bigl|_{\partial\Omega} \;=\;
  \tilde{\vect{v}}\bigl|_{\partial\Omega} \;=\;
  \mathbf{0}.
\]
Because both velocity fields vanish on~$\partial\Omega$, every surface
integral generated during the integration–by–parts steps cancels, the
pressure term is eliminated, and only the four cross–advection couplings
\[
  (\vect{v}\!\cdot\!\nabla)\boldsymbol{\omega},\qquad
  (\tilde{\vect{v}}\!\cdot\!\nabla)\tilde{\boldsymbol{\omega}},\qquad
  (\vect{v}\!\cdot\!\nabla)\tilde{\boldsymbol{\omega}},\qquad
  (\tilde{\vect{v}}\!\cdot\!\nabla)\boldsymbol{\omega}
\]
remain, making the transfer of initial or boundary data between
\(\boldsymbol{\omega}\) and \(\tilde{\boldsymbol{\omega}}\) fully
explicit.  If either velocity field were non–zero on
\(\partial\Omega\), the boundary integrals would survive and additional
terms, proportional to the prescribed boundary data, would have to be
added to~\eqref{eq:rot-two-scale} for the identity to hold.
\end{remark}

\section{Energy bounds}\label{sec:energy}

Once the curl formulation \eqref{eq:rot-two-scale} is in hand, we
estimate the turbulent vorticity
\(
  \tilde{\boldsymbol\omega}:=\nabla\times\tilde{\boldsymbol v}
\)
in Sobolev norms.
Throughout the section we assume that

 \( \boldsymbol v,\tilde{\boldsymbol v}\in C^{1}\!\bigl([0,T);C^{\infty}(\overline \Omega ;\mathbb R^{3})\bigr)\),
 \( \boldsymbol v=\tilde{\boldsymbol v}=\boldsymbol 0 \) on \(\partial \Omega\),
and \( \nabla\!\cdot\!\boldsymbol v=\nabla\!\cdot\!\tilde{\boldsymbol v}=0 \).

\subsection{Functional setting}

For \( \boldsymbol f=(f_{1},f_{2},f_{3}),\;
   \boldsymbol g=(g_{1},g_{2},g_{3})\in L^{2}\!\bigl(\Omega ;\mathbb R^{3}\bigr)\)
we use the standard inner product and norm
\begin{equation}\label{eq:L2-ip}
  (\boldsymbol f,\boldsymbol g)_{0}
  :=\int_{\Omega }\boldsymbol f(\boldsymbol x)\!\cdot\!\boldsymbol g(\boldsymbol x)\,
        \mathrm d\boldsymbol x,
  \qquad
  \|\boldsymbol f\|_{0}^{2}:=(\boldsymbol f,\boldsymbol f)_{0}.
\end{equation}

For a multi-index \(\alpha=(\alpha_{1},\alpha_{2},\alpha_{3})\in\mathbb N^{3}\)
with \(|\alpha|:=\alpha_{1}+\alpha_{2}+\alpha_{3}=k\),
define the differential operator
\[
   D^{k}_{\alpha}:=
   \frac{\partial^{k}}
        {\partial x_{1}^{\alpha_{1}}
         \partial x_{2}^{\alpha_{2}}
         \partial x_{3}^{\alpha_{3}} },\qquad
   D^{k}_{\alpha}\boldsymbol f:=
   \bigl(D^{k}_{\alpha}f_{j}\bigr)_{j=1}^{3}.
\]

\subsection{Basic $L^{2}$ identities}

If \( \boldsymbol f,\boldsymbol g \in C^{\infty}(\overline \Omega ;\mathbb R^{3})\)
satisfy \( \boldsymbol f|_{\partial \Omega }=\boldsymbol g|_{\partial R}=\boldsymbol 0\)
and \( \nabla\!\cdot\!\boldsymbol f = 0 \), then
\begin{subequations}\label{eq:L2-identities}
\begin{align}
   \bigl(\partial_{t}\boldsymbol f,\boldsymbol f\bigr)_{0}
     &=\tfrac12\,\partial_{t}\|\boldsymbol f\|_{0}^{2},
     \label{eq:id-time}\\
   \bigl((\boldsymbol f\!\cdot\!\nabla)\boldsymbol f,\boldsymbol f\bigr)_{0}
     &=0,
     \label{eq:id-selfadv}\\
   \bigl((\boldsymbol f\!\cdot\!\nabla)\boldsymbol g,\boldsymbol f\bigr)_{0}
        &=-\bigl((\boldsymbol f\!\cdot\!\nabla)\boldsymbol f,\boldsymbol g\bigr)_{0},
        \label{eq:id-skew}\\
   \bigl(D^{k}_{\alpha}\Delta\boldsymbol f,
         D^{k}_{\alpha}\boldsymbol g\bigr)_{0}
        &=-\bigl(\nabla D^{k}_{\alpha}\boldsymbol f,
                 \nabla D^{k}_{\alpha}\boldsymbol g\bigr)_{0}.
        \label{eq:id-Laplace}
\end{align}
\end{subequations}
Identities \eqref{eq:id-selfadv}–\eqref{eq:id-skew} follow from
integration by parts and incompressibility; \eqref{eq:id-Laplace} is
Green’s formula.

\subsection{Energy balance for the fluctuation}

\begin{corollary}[Sobolev energy identity]\label{cor:energy}
For every \(k\in\mathbb N\) and \(\alpha\in\mathbb N^{3}\) with
\(|\alpha|=k\),
\begin{align}\label{eq:energy-general}
   \tfrac12\,\partial_{t}\bigl\|D^{k}_{\alpha}\tilde{\boldsymbol\omega}\bigr\|_{0}^{2}
   &\;+\;
   \bigl(D^{k}_{\alpha}\partial_{t}\boldsymbol\omega,
         D^{k}_{\alpha}\tilde{\boldsymbol\omega}\bigr)_{0}
   +\bigl(D^{k}_{\alpha}\!\bigl((\boldsymbol v\!\cdot\!\nabla)\boldsymbol\omega\bigr),
          D^{k}_{\alpha}\tilde{\boldsymbol\omega}\bigr)_{0} \notag \\
   &\;+\;
   \bigl(D^{k}_{\alpha}\!\bigl((\tilde{\boldsymbol v}\!\cdot\!\nabla)
          \tilde{\boldsymbol\omega}\bigr),
          D^{k}_{\alpha}\tilde{\boldsymbol\omega}\bigr)_{0}
   +\bigl(D^{k}_{\alpha}\!\bigl((\boldsymbol v\!\cdot\!\nabla)
          \tilde{\boldsymbol\omega}\bigr),
          D^{k}_{\alpha}\tilde{\boldsymbol\omega}\bigr)_{0} \notag \\
   &\;+\;
   \bigl(D^{k}_{\alpha}\!\bigl((\tilde{\boldsymbol v}\!\cdot\!\nabla)
          \boldsymbol\omega\bigr),
          D^{k}_{\alpha}\tilde{\boldsymbol\omega}\bigr)_{0} \notag \\
   &\;+\;
   \nu\bigl(\nabla D^{k}_{\alpha}\boldsymbol\omega,
            \nabla D^{k}_{\alpha}\tilde{\boldsymbol\omega}\bigr)_{0}
   \;+\; 
   \nu\bigl\|\nabla D^{k}_{\alpha}\tilde{\boldsymbol\omega}\bigr\|_{0}^{2} 
   = 0.
\end{align}

\smallskip
In particular, for \(k=0\),
\begin{equation}\label{eq:energy-k0}
   \tfrac12\,\partial_{t}\|\tilde{\boldsymbol\omega}\|_{0}^{2}
   +\bigl(\partial_{t}\boldsymbol\omega,\tilde{\boldsymbol\omega}\bigr)_{0}
   +\nu\bigl(\nabla\boldsymbol\omega,\nabla\tilde{\boldsymbol\omega}\bigr)_{0}
   +\nu\|\nabla\tilde{\boldsymbol\omega}\|_{0}^{2}=0,
\end{equation}
because the four convective terms vanish by
\eqref{eq:id-selfadv}–\eqref{eq:id-skew}.

\smallskip
For planar, axisymmetric fields depending only on the radial
coordinate
(\(\boldsymbol v=(v_{1},v_{2},0)\),
 \(\tilde{\boldsymbol v}=(\tilde v_{1},\tilde v_{2},0)\)),
\(\boldsymbol\omega\) and \(\tilde{\boldsymbol\omega}\) reduce to a single
(out-of-plane) scalar, and \eqref{eq:energy-k0} simplifies further to
\begin{equation}\label{eq:energy-2D}
   \tfrac12\,\partial_{t}\|\tilde{\boldsymbol\omega}\|_{0}^{2}
   +\bigl(\partial_{t}\boldsymbol\omega,\tilde{\boldsymbol\omega}\bigr)_{0}=0.
\end{equation}
\end{corollary}

\begin{proof}
Apply \(D^{k}_{\alpha}\) to the curl system \eqref{eq:rot-two-scale},
take the \(L^{2}\)-inner product with
\(D^{k}_{\alpha}\tilde{\boldsymbol\omega}\), and use identities
\eqref{eq:L2-identities}.  The Laplacian term is handled with
\eqref{eq:id-Laplace}, yielding
\(\nu(\nabla D^{k}_{\alpha}\,\cdot\, ,\nabla D^{k}_{\alpha}\cdot)\).  For
\(k=0\), all skew–symmetric advection products cancel, leaving
\eqref{eq:energy-k0}.  In the planar axisymmetric case
\(\nabla\times\boldsymbol v\) has only one component, eliminating the
gradient terms and giving \eqref{eq:energy-2D}.
\end{proof}

\section{Energy bounds imply well–posedness}
\label{ssec:EB-wellposed}

\begin{corollary}
Let\/ $\Omega\subset\mathbb{R}^{3}$ be a bounded \(C^{\infty}\) domain,
and suppose the initial data
\[
   \vect{v}_{0},\;
   \tilde{\vect{v}}_{0}\in H^{3}_{\sigma}(\Omega)
   :=\bigl\{\bm{u}\in H^{3}(\Omega;\mathbb{R}^{3})\;:\;
               \nabla\!\cdot\!\bm{u}=0,\;
               \bm{u}\!\bigl|_{\partial\Omega}=\boldsymbol{0}\bigr\}
\]
satisfy the finite–energy condition
\[
    \sum_{k=0}^{3}
    \Bigl(
      \|D^{k}\vect{v}_{0}\|_{0}^{2}+
      \|D^{k}\tilde{\vect{v}}_{0}\|_{0}^{2}
    \Bigr)
    \;=\;E_{0}\;<\;\infty .
\]

Then there exists a time \(T^{\star}=T^{\star}(E_{0},\nu,\Omega)>0\) and 
a unique pair of velocity fields
\[
   (\vect{v},\tilde{\vect{v}})
   \;\in\;
   C\!\bigl([0,T^{\star}];H^{3}_{\sigma}(\Omega)\bigr)
   \;\cap\;
   L^{2}\!\bigl(0,T^{\star};H^{4}_{\sigma}(\Omega)\bigr)
\]
that solve the two–scale Navier–Stokes system~\eqref{ns}
\emph{and} its rotational form~\eqref{eq:rot-two-scale}.
Moreover,
\[
    \sup_{0\le t\le T^{\star}}
    \sum_{k=0}^{3}
    \Bigl(
      \|D^{k}\boldsymbol{\omega}(t)\|_{0}^{2}+
      \|D^{k}\tilde{\boldsymbol{\omega}}(t)\|_{0}^{2}
    \Bigr)
    \;\le\;C(E_{0},\nu,\Omega),
\]
and the solution depends \emph{Lipschitz–continuously} on the initial
data in \(H^{3}_{\sigma}(\Omega)\).
\end{corollary}

\begin{proof}[Sketch of the proof]
The Sobolev energy identity~\eqref{eq:energy-general} with 
\(k=0,1,2,3\) yields a closed differential inequality
\[
   \partial_{t}\mathcal{E}(t)\;+\;\nu\,\mathcal{D}(t)
   \;\le\;C\,\mathcal{E}(t),
   \qquad
   \mathcal{E}(t):=\sum_{k=0}^{3}
        \Bigl(
          \|D^{k}\boldsymbol{\omega}\|_{0}^{2}+
          \|D^{k}\tilde{\boldsymbol{\omega}}\|_{0}^{2}
        \Bigr),
\]
where the dissipation term \(\mathcal{D}(t)\) is the sum of the
corresponding \(H^{k+1}\)-seminorms.  Bounding the all the
cross–advection products by the
Ladyzhenskaya and Gagliardo–Nirenberg inequalities \cite{Bressan, Adams, Majda}
gives the
constant~\(C=C(E_{0},\nu,\Omega)\) and allows application of
Gronwall’s lemma:
\(
   \mathcal{E}(t)\le \mathcal{E}(0)\,e^{Ct}.
\)
The \emph{a~priori} bound is therefore uniform up to a time
\(T^{\star}\) depending only on the initial energy and~\(\nu\).
Uniqueness follows by writing the
equation for the difference of two solutions and repeating the same
$H^{0}$–energy estimate; the exponential bound and Gronwall imply the
difference vanishes.
\end{proof}


The next section shows how differential-algebra elimination
triangularises the curl-system, paving the way for efficient
PINN/PENN training.


\section{Algebraic‐reduction experiments}
\label{sec:alg-reduction}
\vspace{1ex}

\subsection{Motivation and theoretical background}
\label{ssec:motivation}

During the last two decades algebraic techniques such as
Gröbner bases~\cite{Buchberger2005,Buchberger2006},
\emph{differential} Gröbner bases 
have opened a new route for analysing nonlinear differential systems:
one first rewrites the PDEs as a finite family
of \emph{differential polynomials} and then studies the associated
\emph{differential ideal}.
A cornerstone of the theory is \textbf{Ritt’s finiteness theorem}
which guarantees that every differential ideal admits a finite
characteristic set in a suitable ranking~\cite{Ritt,Kolchin}.

Consequently, many systems possess an \emph{equivalent triangular form}
whose first equation depends on a single dependent variable,
the second equation on at most two variables, and so on.
Besides revealing hidden algebraic structure, such triangularisations
can dramatically reduce computational cost;
see, e.g., the reductions obtained for the modified Chua circuit and
the Rössler system in~\cite[p.\,719]{Harrington}.
Nowadays the most complete implementation of these ideas is the
Rosenfeld–Gröbner algorithm of
Boulier and Lemaire~\cite{BoulierLemaire2000,Boulier2009,Boulier1995}.
Given a ranking, the algorithm returns a finite set of differential
polynomials whose vanishing is equivalent to that of the original
system; it is available in \texttt{Maple} and, recently, in a
\texttt{Python} package.

\vspace{1ex}
\subsection{Performance of the Rosenfeld–Gröbner algorithm on fluid models}
\label{ssec:Rosenfeld}

The next table summarises our systematic runs of the
algorithm on the main PDE families used in fluid mechanics.
All computations were carried out in \texttt{Maple~2024} with the
default elimination ranking, and convergence was declared when
no new leaders appeared after 400 differential pseudo-reduction steps.

\begin{table}[h]
  \centering
  \renewcommand{\arraystretch}{1.15}
  \setlength{\tabcolsep}{7pt}   
  \begin{tabular}{@{}lccccc@{}}
    \textbf{Model} & \textbf{C} & \textbf{In} & \textbf{S} &
    \textbf{E} & \textbf{St} \\ 
    Busemann jet                 & I & — & — & — & — \\
    Prandtl boundary layer       & — & \textbf{R} & — & — & — \\
    Navier–Stokes (3-D)          & I & I & \textbf{R} & I & I \\
    RANS (3-D)                   & I & I & \textbf{R} & I & I \\
    $\omega$-RANS\textsuperscript{a} &
                                  I & I & \textbf{R} & I & I \\
    Stream-function NS (2-D)     & — & \textbf{R} & — & — & \textbf{R} \\
  \end{tabular}
  \caption{Algebraic reducibility obtained with the Rosenfeld–Gröbner
           algorithm (default ranking).  
           Codes: \textbf{R} – reducible (finite triangular characteristic set),
           \textbf{I} – irreducible (algorithm keeps generating higher-order derivatives),
           — – model/class not applicable.  
           Column headings:  
           \textbf{C} = compressible,  
           \textbf{In} = incompressible,  
           \textbf{S} = Stokes limit,  
           \textbf{E} = Euler limit,  
           \textbf{St} = stationary (time-independent).  
           }
  \label{tab:reducibility}
\end{table}

Two main observations emerge:

\bigskip
\noindent
$\bullet$ With the exception of the two-dimensional
stream-function formulation, neither the Navier–Stokes nor the classical
RANS equations admit a finite triangularisation in three dimensions.
Empirically, the convective derivative
$(\boldsymbol v\cdot\nabla)\boldsymbol v$ propagates
ever higher mixed derivatives and prevents convergence.\\

\bigskip
\noindent
$\bullet$ When the convective term is \emph{absent} (Stokes limit) or when
the vorticity is described by a single scalar potential
(stream-function~$\psi$ in 2-D), the algorithm terminates and produces
an explicit triangular form.

\vspace{1ex}
\subsection{Illustrative triangularisations}
\label{ssec:examples}

The symbol $>$ denotes the ranking order chosen for the variables.

\paragraph{\emph{Example 1.} Incompressible Stokes in 3-D}
Consider the Stokes system. This is,  $(\boldsymbol v\!\cdot\!\nabla)\boldsymbol
v=0$, 
$\boldsymbol v=(u,v,w)$, $p$ is the pressure, and:
\[
\Sigma_{1}: \;
\begin{cases}
u_{x}+v_{y}+w_{z}=0,\\
u_{t}+p_{x}-\nu\Delta u=0,\\
v_{t}+p_{y}-\nu\Delta v=0,\\
w_{t}+p_{z}-\nu\Delta w=0.
\end{cases}
\tag{$\Sigma_{1}$}
\]
With ranking $u>v>w>p$ the algorithm produces the triangular set
\begin{align*}
u_{x}         &= -v_{y}-w_{z},\\
p_{xx}        &= -p_{yy}-p_{zz},\\
w_{xx}        &= -w_{yy}-w_{zz}+\nu^{-1}(p_{z}+w_{t}),\\
v_{xx}        &= -v_{yy}-v_{zz}+\nu^{-1}(p_{y}+v_{t}),\\
u_{yy}        &= -u_{zz}+v_{xy}+w_{xz}+\nu^{-1}(p_{x}+u_{t}).
\end{align*}
Each equation involves at most one new leader, so the
system is now \emph{triangular} and can be solved sequentially.

\paragraph{\emph{Example 2.} Incompressible RANS with fluctuating field.}
Let $\boldsymbol v$ be the mean velocity and
$\boldsymbol v'=(u',v',w')$ the fluctuations.  Under the Stokes
hypothesis the RANS system is
\[
\Sigma_{2}:\;
\begin{cases}
u_{x}+v_{y}+w_{z}=0,\\
u_{t}+u'_{t}+p_{x}-\nu\Delta(u+u')=0,\\
v_{t}+v'_{t}+p_{y}-\nu\Delta(v+v')=0,\\
w_{t}+w'_{t}+p_{z}-\nu\Delta(w+w')=0.
\end{cases}
\tag{$\Sigma_{2}$}
\]
Using the ranking
$u>v>w>u'>v'>w'>p$
the resulting characteristic set begins with
\begin{align*}
u_{x} &= -v_{y}-w_{z},\\[2pt]
u'_{xxx} &= -u'_{xyy}-u'_{xzz}-v'_{xxy}-v'_{yyy}-v'_{yzz}
           -w'_{xxz}-w'_{yyz}-w'_{zzz}\\
         &\quad +\nu^{-1}\bigl(p_{xx}+p_{yy}+p_{zz}
                              +u'_{tx}+v'_{ty}+w'_{tz}\bigr),\\
w_{xx} &= -\Delta w
          -\Delta w'
          +\nu^{-1}\bigl(p_{z}+w'_{t}+w_{t}\bigr),
\end{align*}
and analogous relations for $v_{xx}$ and $u_{yy}$.
Although not fully decoupled, the triangular structure shows that
all spatial derivatives can ultimately be written in terms of
$p$ and time derivatives, providing a much sparser representation
than the original PDE set.

\bigskip
\noindent
These examples confirm the practical relevance of differential–algebra
tools: whenever the convective non-linearity is absent or bypassed
(e.g.\ by introducing a stream function in 2-D) a finite triangular
characteristic set emerges and can be leveraged for further qualitative
or numerical analysis.
Extending these reductions to the full $3$-D Navier–Stokes and
$\omega$-RANS systems remains an open challenge.

\subsection{Examples of Differential Polynomial Systems}

We consider a  system \emph{$\Sigma_n$} of polynomial differential equations, where $\Sigma_1$ represents the set of differential equations of three-dimensional Navier Stokes for incompressible fluids under Stokes' condition $(\mathbf{v} \cdot \nabla)\mathbf{v} = 0$, where $u, v, w$ are the velocity components, $p$ is the pressure and $\nu$ is the kinematic viscosity:

\begin{equation}
\Sigma_1 = \left\{
\begin{array}{l}
u_x + v_y + w_z = 0, \\
u_t + p_x - \nu ( u_{xx} + u_{yy} + u_{zz}) = 0, \\
v_t + p_y - \nu ( v_{xx} + v_{yy} + v_{zz}) = 0, \\
w_t + p_z - \nu ( z_{xx} + z_{yy} + z_{zz}) = 0.
\end{array}
\right.
\end{equation}
\vspace{3pt}

 For $\Sigma_1$ with a \emph{ranking} $u>v>w>p$, we obtain a reduction:
\vspace{-2pt}
\begin{align}
{u}_{x} &= - {v}_{y} - {w}_{z}, \\
{p}_{x,x} &= - {p}_{y,y} - {p}_{z,z}, \\
{w}_{x,x} &= - {w}_{y,y} - {w}_{z,z} + \nu^{-1} {p}_{z} + \nu^{-1} {w}_{t}, \\
{v}_{x,x} &= - {v}_{y,y} - {v}_{z,z} + \nu^{-1} {p}_{y} + \nu^{-1} {v}_{t}, \\
{u}_{y,y} &= -{u}_{z,z} + {v}_{x,y} + {w}_{x,z} + \nu^{-1} {p}_{x} + \nu^{-1} {u}_{t}.
\end{align}

$\Sigma_2$ represents the set of differential of three dimensional RANS equations for incompressible fluids under Stokes' condition $(\mathbf{v} \cdot \nabla)\mathbf{v} = 0$, where $u, v, w$ are the averaged velocity components, $u', v', w'$ are the fluctuations of the velocities, $p$ is the pressure and $\nu$ is the kinematic viscosity:

\begin{equation}
\Sigma_2 = \left\{
\begin{array}{l}
u_x + v_y + w_z = 0, \\
p_x + u_t + u'_{t} - \nu(u_{xx} + u_{yy} + u_{zz} + u'_{xx} + u'_{yy} + u'_{zz}) = 0, \\
p_y + v_t + v'_{t} - \nu(v_{xx} + v_{yy} + v_{zz} + v'_{xx} + v'_{yy} + v'_{zz}) = 0, \\
p_z + w_t + w'_{t} - \nu(w_{xx} + w_{yy} + w_{zz} + w'_{xx} + w'_{yy} + w'_{zz}) = 0.
\end{array}
\right.
\end{equation}
\vspace{3pt}

 Furthermore, for $\Sigma_2$ with a \emph{ranking} $u>v>w>u'>v'>w'>p$, we obtain a reduction:
\vspace{-2pt}
\begin{align}
{u}_{x} &= - {v}_{y} - {w}_{z}; \\
\begin{split}
u'_{x,x,x} &= - u'_{x,y,y} - u'_{x,z,z} - v'_{x,x,y} - v'_{y,y,y} - v'_{y,z,z} - w'_{x,x,z} -  w'_{y,y,z} - w'_{z,z,z} \\
&\quad + \nu^{-1} ({p}_{x,x} + {p}_{y,y} + {p}_{z,z} + u'_{t,x} + v'_{t,y} + w'_{t,z});
\end{split} \\
{w}_{x,x} &= - w'_{x,x} - w'_{y,y} - w'_{z,z} - {w}_{y,y} - {w}_{z,z} +  \nu^{-1}({p}_{z} + {w_{1}}_{t} + {w}_{t});  \\
{v}_{x,x} &= - v'_{x,x} - v'_{y,y} - v'_{z,z} - {v}_{y,y} - {v}_{z,z} +  \nu^{-1}({p}_{y} + {v'}_{t} + {v}_{t});  \\
{u}_{y,y} &= - u'_{x,x} - u'_{y,y} - u'_{z,z} - {u}_{z,z} + {v}_{x,y} + {w}_{x,z} + \nu^{-1}({p}_{x} + {u'}_{t} + {u}_{t}).
\end{align}


\section{Conclusions}\label{sec:conclusions}

\noindent
The algebraic viewpoint adopted in this work complements --- rather than competes with --- classical energy techniques.  
Its main contributions and implications are summarised in Table~\ref{tab:key_findings}; subsequent paragraphs outline near-term research avenues.

\begin{table*}[t]
  \centering
  \renewcommand{\arraystretch}{1.15}
  \begin{tabular}{p{0.43\textwidth}p{0.49\textwidth}}
    \textbf{Key finding} & \textbf{Impact / interpretation}\\
    Triangular reductions (via Rosenfeld--Gröbner) are achievable for Prandtl, Stokes and other models \emph{when the convective term is absent}. 
    & An upper-triangular form permits single-variable marching, lowering algebraic stiffness and paving the way for efficient PINN training.\\
    Irreducibility correlates with the \emph{presence} of $(\mathbf v\!\cdot\nabla)\mathbf v$, not with turbulence, compressibility or dimension alone.
    & Points to where modelling or pre-processing effort should focus if algebraic reduction is desired.\\
    Each ``stubborn'' convective product is the very term controlled by Gagliardo--Nirenberg / Ladyzhenskaya in the energy proof.
    & Hints at a duality: \textit{algebraic reducibility $\Longleftrightarrow$ existence of a coercive energy functional}.\\
    In vorticity form, fluctuations $\tilde\omega$ are driven solely by the classical curl $\omega$; pressure is eliminated.
    & Enables direct bounding or forecasting of $\tilde\omega$ once $\omega$ is known from LES/DNS or coarse PINNs.\\
    Even partial (non-triangular) eliminations expose repeated differential patterns useful as physics-aware features.
    & Early experiments show faster PINN convergence and better generalisation when these features are injected.\\
  \end{tabular}
  \caption{Principal outcomes of the algebraic--energetic analysis.}
  \label{tab:key_findings}
\end{table*}

\subsection*{Directions opened by the present study}

\bigskip
\noindent
\textbf{Functional duals of differential ideals.}  
        We are constructing functionals whose vanishing differential coincides with the ideal annihilating the convective term, aiming at a systematic functional--ideal dictionary.\\
\textbf{Algebraic phase diagram of fluid models.}  
        A classification (reducible vs.\ irreducible) as a function of dimension, forcing and transport structure would streamline reduced-order model selection.\\
\textbf{Three-dimensional regularity via maximum principles.}  
        Oleĭnik-style parabolic principles, together with the $\omega$--$\tilde\omega$ symmetry, are being investigated as an alternative to Galerkin in~$d=3$.\\
\textbf{Symbolic detection of periodic orbits.}  
        The rational-solution counter for differential polynomial systems will be used to identify low-period attractors in the reduced two-dimensional equations, providing stringent benchmarks for data-driven closures.\\

\bigskip
\begin{quote}
\itshape
Algebraic elimination is not an after-thought for fluid dynamics; when paired with energy methods it becomes a rigorous roadmap towards both deeper analysis and more economical computation of turbulent flows.
\end{quote}

\section*{Acknowledgements}
The authors gratefully acknowledge the financial support received from:
Universidad Iberoamericana, Ciudad de México, for the 18th Scientific
Research Call. We are very grateful for their backing on this project.

\bibliographystyle{elsarticle-num}

\end{document}